# The degenerate special Lagrangian equation on Riemannian manifolds

Matthew Dellatorre

**Abstract**

We show that the degenerate special Lagrangian equation, recently introduced by Rubinstein–Solomon, induces a global equation on every Riemannian manifold, and that for certain associated geometries this equation governs, as it does in the Euclidean setting, geodesics in the space of positive Lagrangians. For example, geodesics in the space of positive Lagrangian sections of a smooth Calabi–Yau torus fibration are governed by the Riemannian DSL on the base manifold.

We then develop their analytic techniques, specifically modifications of the Dirichlet duality theory of Harvey–Lawson, in the Riemannian setting to obtain continuous solutions to the Dirichlet problem for the Riemannian DSL and hence continuous geodesics in the space of positive Lagrangians.

## 1 INTRODUCTION

Let $f \in C^2([0,1] \times \mathbb{R}^n)$ and $\theta \in (-\pi, \pi]$. Then $f(t,x)$ satisfies the *degenerate special Lagrangian equation of phase* $\theta$ if

$$\operatorname{Im}\left(e^{-\sqrt{-1}\theta}\det(I_n + \sqrt{-1}\nabla^2 f)\right) = 0 \quad \text{and} \quad \operatorname{Re}\left(e^{-\sqrt{-1}\theta}\det(I + \sqrt{-1}\nabla_x^2 f)\right) > 0. \tag{1}$$

Here $I_n$ denotes the diagonal $(n+1) \times (n+1)$ matrix with diagonal entries $(0, 1, \ldots, 1)$. The degenerate special Lagrangian equation (DSL) was introduced by Rubinstein–Solomon [RS] in connection to geodesics in the space of positive Lagrangians of a Calabi–Yau manifold. It is a fully nonlinear, degenerate elliptic equation.

When the featured Calabi–Yau is $\mathbb{C}^n$, the geodesic endpoint problem in the space of positive graph Lagrangians corresponds to solving the Dirichlet problem for the DSL. In particular, the conditions in (1) capture, respectively, the notions of geodesic and positivity in this setting. Under appropriate boundary conditions, unique continuous solutions to the Dirichlet problem for the DSL exist. This was accomplished in [RS] by finding a natural notion of subsolution to the DSL and then adapting the Dirichlet duality framework of Harvey–Lawson [HL2] for degenerate elliptic equations in Euclidean space.

Harvey–Lawson [HL3] have also developed a Dirichlet duality theory for equations on Riemanian manifolds. The starting point for this framework is an equation $F$ in Euclidean space and a Riemannian manifold $M$. Assuming the topology on $M$ is sufficiently mild and the symmetry of $F$ is sufficiently high, one can define a global equation on $M$ that is locally modelled on $F$. Thus, from this point of view, it is natural to consider the equation induced by the DSL on Riemannian manifolds, and that is the purpose of this note.

Geodesics in the space of positive Lagrangians play a crucial role in a program initiated by Solomon [S1, S2] (see also [SY]) to understand the existence and uniqueness of special Lagrangian submanifolds in Calabi–Yau manifolds. They also play a key role in a new approach to the Arnold conjecture put forth by Rubinstein–Solomon [RS, Section 2.3] and in the development of a pluripotential theory for Lagrangian graphs initiated in [RS] (see also [DR]).

### 1.1 Results

We begin with the easy observation that the DSL carries over to a global equation on every Riemannian manifold. For precise statements see Section 2 and Appendix A.

**Proposition 1.1.1.** *Given any Riemannian manifold $M$, the degenerate special Lagrangian equation carries over (in the sense of Harvey–Lawson) to a global equation on $\mathbb{R} \times M$, locally modelled on the DSL. We refer to this equation as the Riemannian DSL on $\mathbb{R} \times M$.*

The second goal of this paper is to show that the Riemannian DSL also governs geodesics in the space of positive Lagrangians for certain geometries associated to $M$; for example when the ambient Calabi–Yau manifold is the cotangent bundle $T^*M$ of certain paralellizable manifolds or when it is a Lagrangian fibration over $M$. See Section 3 for details and definitions.

**Theorem 1.1.2.** *Let $M$ be integrably parallelizable. Then $T^*M$ admits a Calabi–Yau structure, and the Riemannian DSL on $[0,1] \times M$ governs geodesics in the space of positive graph Lagrangians in $T^*M$.*

**Theorem 1.1.3.** *Let $X$ be a smooth Calabi–Yau torus fibration over $M$. Then the Riemannian DSL on $[0,1] \times M$ governs geodesics in the space of positive Lagrangian sections in $X$.*

This motivates the third goal of this paper which is to solve the Dirchlet problem for the Riemannian DSL, as this corresponds to the endpoint problem for geodesics. More specifically, we aim to solve the Dirichlet problem on domains of the form

$$\mathcal{D} = [0,1] \times D \subset \mathbb{R} \times M,$$

where $D \subset M$ is a bounded domain. This is accomplished by following the approach of Rubinstein–Solomon. In particular, we extend the Dirichlet duality theory of Harvey–Lawson to include certain domains with corners in Riemannian manifolds, such as $\mathcal{D}$ when $D \subset M$ has boundary. This extension is contained in Theorem 4.2.1.

Under appropriate boundary conditions, Theorem 4.2.1 provides continuous solutions to the Dirichlet problem for the Riemannian DSL on $\mathcal{D}$ and thus continuous geodesics in the space of positive Lagrangians. A special case of this result is the following theorem. By strictly convex we mean that all of the eigenvalues of the second fundamental form $II_{\partial D}$ are strictly positive, and by admissible we mean that this local frame for $TM$ is part of a family of frames whose transition maps are $O_n$-valued (see Appendix A).

**Theorem 1.1.4.** *Let $(M,g)$ be a complete simply-connected Riemannian manifold with non-positive sectional curvature, and let $D \subset M$ be a bounded strictly convex domain. For $i = 0, 1$, let $\phi_i \in C^2(D)$ satisfy*

$$\operatorname{tr} \tan^{-1} (\operatorname{Hess} \phi_i(e,e)) \in (c - \pi/2, c + \pi/2), \tag{2}$$

*where $e = (e_1, \cdots, e_n)$ is any admissible local frame for the tangent bundle. Then there exists a unique solution $u \in C^0(\overline{\mathcal{D}})$ to the Dirichlet problem for the Riemannian DSL of phase $\theta$ (where $c = \theta \mod 2\pi$) with $u|_{\{i\} \times D} = \phi_i$ and $u|_{[0,1] \times \partial D}$ affine in $t$.*

See Remark 5.0.2 for a proof. In certain settings (see Section 3 and [RS, Section 1]), condition (2) is equivalent to the condition that the graph of $d\phi_i$ in the cotangent bundle is a positive Lagrangian.

## 1.2 Geometry of the space of Lagrangians

The following section is based on the work of Solomon [S1, S2] and briefly recalls the terminology concerning the geometry of the space of positive Lagrangians.

Let $L$ be an $n$-dimensional real manifold and $(X, J, \omega, \Omega)$ an almost Calabi–Yau manifold of complex dimension $n$. That is, $(X, J, \omega)$ is a Kähler manifold and $\Omega$ is a nowhere vanishing holomorphic $n$-form. Define

$$\mathcal{L} = \{\Gamma \subset X : \Gamma \text{ is an oriented Lagrangian submanifold diffeomorphic to } L\}.$$

For $\theta \in (-\pi, \pi]$, the space of $\theta$-positive Lagrangians is defined as

$$\mathcal{L}_\theta^+ = \{\Gamma \in \mathcal{L} \mid \operatorname{Re}(e^{-\sqrt{-1}\theta}\Omega)|_\Gamma > 0\}. \tag{3}$$

Denote by $\mathcal{O}_\theta \subset \mathcal{L}_\theta^+$ a connected component of the intersection of $\mathcal{L}_\theta^+$ with an orbit $\operatorname{Ham}(X,\omega)$ acting on $\mathcal{L}$, where $\operatorname{Ham}(X,\omega)$ is the group of compactly supported Hamiltonion diffeomorphisms of $X$.

When $L$ is compact the tangent space to $\mathcal{O}_\theta$ can be identified with the space of smooth functions satisfying a normalization condition

$$T_\Gamma \mathcal{O}_\theta \equiv \{h \in C^\infty(\Gamma) \mid \int_\Gamma h \operatorname{Re} \Omega = 0\}, \tag{4}$$

and a weak Riemannian metric on $\mathcal{O}_\theta$ is defined by

$$(h, k)_\theta|_\Gamma := \int_\Gamma hk \operatorname{Re}(e^{-\sqrt{-1}\theta}\Omega|_\Gamma), \text{ for } h, k \in T_\Gamma \mathcal{O}_\theta. \tag{5}$$

When $L$ is non-compact the normalization condition in (4) can be dropped and the tangent space at $\Gamma$ is isomorphic to the space of compactly supported functions on $\Gamma$.

More specifically, given a path $\Lambda : [0,1] \to \mathcal{O}_\theta$ and a family of diffeomorphisms $g_t : L \to \Lambda_t$, let $h_t : \Lambda_t \to \mathbb{R}$ be the unique function satisfying

$$g_t^* \iota_{\frac{dg_t}{dt}} \omega = d(h_t \circ g_t), \tag{6}$$

and the normalization condition in (4). Then the velocity vector to $\Lambda$ is defined as $\frac{d\Lambda_t}{dt} \equiv h_t$.

Given a vector field $q_t \in T_{\Lambda_t}\mathcal{O}_\theta$ along $\Lambda$, the Levi–Civita coavariant derivative of $q_t$ in the direction of $\frac{d\Lambda_t}{dt}$ is defined by

$$\frac{Dq_t}{dt} = \left( \frac{\partial}{\partial t}(q_t \circ g_t) + g_t^* dq_t(\zeta_t) \right) \circ g_t^{-1}, \tag{7}$$

where $\zeta_t$ is the unique vector field on $L$ such that

$$\iota_{\zeta_t} g_t^* \operatorname{Re}(e^{-i\theta}\Omega) = -g_t^* \iota_{\frac{dg_t}{dt}} \operatorname{Re}(e^{-i\theta}\Omega), \tag{8}$$

viewing $g_t : L \to \Gamma_t \subset X$ as a map from $L$ to $X$.

The geodesic equation is then found by taking $q_t = h_t = \frac{d\Lambda_t}{dt}$:

$$\frac{Dh_t}{dt} = \left( \frac{\partial}{\partial t}(h_t \circ g_t) + g_t^* dh_t(\zeta_t) \right) \circ g_t^{-1} = 0. \tag{9}$$

## 1.3 Dirichlet duality theory

This section provides a brief sketch of Harvey–Lawson's [HL2, HL3] Dirichlet duality theory. See Appendix A for a more thorough review. See also the work of Slodkowski [Sl1, Sl2, Sl3, Sl4].

Let $\Omega$ be a bounded domain in $\mathbb{R}^n$ and consider an equation of the form $F(D^2 u) = 0$ on $\Omega$. To any such equation $F$, Dirichlet duality associates a *subequation* $\mathcal{F}$. This is a closed proper subset of $\operatorname{Sym}^2(\mathbb{R}^n)$ that is invariant under translation by positive matrices. Roughly speaking,

$$F(D^2 u) = 0 \longrightarrow \mathcal{F} := \{A \in \operatorname{Sym}^2(\mathbb{R}^n) \ : \ F(A) \geq 0\}.$$

In most cases, the set $\mathcal{F}$ will usually be a more regular proper subset of this set. Regardless, solutions $u \in C^2(\Omega)$ of $F$ must satisfy $D_x^2 u \in \partial \mathcal{F}$.

This gives rise to a natural notion of a subsolution to $F$. A $C^2$ function $u$ is $\mathcal{F}$-subharmonic on $\Omega$ if

$$D_x^2 u \in \mathcal{F}, \quad \forall x \in \Omega. \tag{10}$$

This definition extends to upper semi-continuous functions in a viscosity-like way via $C^2$ test functions, and these $\mathcal{F}$-subharmonic functions comprise the subsolutions.

The class of $\mathcal{F}$-subharmonic functions remarkably share most of the important properties that the classical subharmonic and convex functions satisfy. For example, closure under decreasing limits and taking maxima, decreasing limits, uniform limits, and upper envelopes. See Theorem A.3.1.

For each subequation $\mathcal{F}$, there is an associated *dual* subequation $\widetilde{\mathcal{F}}$, defined as $\widetilde{\mathcal{F}} = -(\sim \operatorname{Int}\mathcal{F})$. The importance of the dual subequation lies in the fact that

$$-\widetilde{\mathcal{F}} \cap \mathcal{F} = \partial \mathcal{F} \subset \{F = 0\}. \tag{11}$$

This immediately allows one to define a notion of weak solution. Let $u$ be a $C^2$ function such that $u$ is $\mathcal{F}$-subharmonic and $-u$ is $\widetilde{\mathcal{F}}$-subharmonic on $\Omega$. Then, by (10) and (11), $u$ satisfies

$$F(D^2 u(x)) = 0 \quad \forall x \in \Omega.$$

A function $u$ is said to be $\mathcal{F}$-*harmonic* if $u$ is $\mathcal{F}$-subharmonic and $-u$ is $\widetilde{\mathcal{F}}$-subharmonic. These functions will comprise our weak solutions. Note that since both $u$ and $-u$ are upper semi-continuous, $u$ is automatically continuous.

Given a domain $\Omega$ and an equation $F$, the existence of continuous solutions to Dirichlet problem for $F$ requires that $\partial \Omega$ is $\mathcal{F}$-*convex*, where $\mathcal{F}$ is a subequation associated to $F$. This subequation-specific convexity generalizes conventional convexity. The Perron method is used to construct solutions, and this convexity implies the existence of barrier functions. This explicit relationship between the equation and the necessary geometry of the boundary is one of the exceptional aspects of Dirichlet duality.

## 1.4 Organization

In the next section, we prove Proposition 1.1.1, showing that the DSL induces (in the sense of Harvey–Lawson) a global equation on every Riemannian manifold. We then geometrically motivate our study of the Riemannian DSL in Section 3 by proving Theorem 1.1.2 and Theorem 1.1.3. In Section 5, we extend the Dirichlet duality theory in the Riemannian setting to include domains with corners, proving a generalization of Theorem 1.1.4. In Section 6, we use these results to obtain unique continuous solutions to the Dirichlet problem for the DSL on Riemannian manifolds, and hence continuous geodesics. Finally, for ease of reference, we include an appendix with a brief summary of Dirichlet duality theory.

# 2 THE RIEMANNIAN DSL SUBEQUATION

When $X = \mathbb{C}^n \cong \mathbb{R}^n \oplus \sqrt{-1}\mathbb{R}^n$, with the standard Calabi–Yau structure,

$$\omega = \frac{\sqrt{-1}}{2} \sum_j dz_j \wedge d\bar{z}_j \quad \text{and} \quad \Omega = dz_1 \wedge \ldots dz_n,$$

and $L$ is identified with $\mathbb{R}^n \times \{0\} \subset \mathbb{C}^n$, the analysis in Section 1.2 leads to the degenerate special Lagrangian equation.

**Theorem 2.0.1.** *[RS, Proposition 2.3] Let $\theta \in (-\pi, \pi]$ and let $k_i \in C^2(\mathbb{R}^n)$, $i = 0, 1$ be such that $graph(dk_i) \subset \mathbb{C}^n$ are elements of $\mathcal{O}_\theta$. Let $k \in C^2([0,1] \times \mathbb{R}^n)$ be such that $graph(d_x k(t, \cdot)) \subset \mathbb{C}^n$ is an element of $\mathcal{O}_\theta$ for each $t \in [0,1]$. Then $t \to graph(d_x k(t, \cdot))$ is a geodesic in $(\mathcal{O}_\theta, (\cdot, \cdot))$ with endpoints $graph(dk_i)$, $i = 0, 1$, if and only if $k$ satisfies*

$$\text{Im}\left(e^{-\sqrt{-1}\theta}\det(I_n + \sqrt{-1}\nabla^2 k)\right) = 0 \quad \text{and} \quad \text{Re}\left(e^{-\sqrt{-1}\theta}\det(I + \sqrt{-1}\nabla_x^2 k)\right) > 0, \qquad (12)$$

*and $k(0, \cdot) = k_0 + c$, $k(1, \cdot) = k_1 + c$, for a constant $c \in \mathbb{R}$.*

## 2.1 The DSL subequation

In order to obtain a subequation for the DSL, Rubinstein–Solomon associate to each $u \in C^2(\mathcal{D})$, where $\mathcal{D} = (0, 1) \times D$, the circle valued function

$$\Theta_u(t, x) = \Theta(\nabla^2 u(t, x)) = \arg \det(I_n + \sqrt{-1}\nabla^2 u(t, x)) \in S^1,$$

defined where $\det(I_n + \sqrt{-1}\nabla^2 u(t, x)) \neq 0$. This angle $\Theta$ is called the *space-time Lagrangian angle* by analogy with the Lagrangian angle of Harvey–Lawson [HL1]. Accordingly, if $u \in C^2(\mathcal{D})$ solves the DSL of phase $\theta$, then $\Theta_u \equiv \theta$.

For a complex matrix $B$, let $\text{spec}(B)$ be its set of eigenvalues, and for $\lambda \in \text{spec}(B)$ denote by $m(\lambda)$ its multiplicity as a root of the characteristic polynomial. Let $\mathcal{S} \subset \text{Sym}^2(\mathbb{R}^{n+1})$ denote the set of symmetric matrices with all zeros in the first row and column, and for $A \in \text{Sym}^2(\mathbb{R}^{n+1}) \setminus \mathcal{S}$, define

$$\widehat{\Theta}(A) = \sum_{\lambda \in \text{spec}(I_n + \sqrt{-1}A)} m(\lambda) \arg(\lambda),$$

branch of arg with values in $(-\pi, \pi]$. Note that $\arg \det(I_n + \sqrt{-1}A) = \widehat{\Theta}(A) \mod 2\pi$. Denote by $\widetilde{\Theta}$ the minimal upper semi-continuous extension of $\widehat{\Theta}$ to $Sym^2(\mathbb{R}^{n+1})$.

**Theorem 2.1.1.** *[RS, Theorem 1.1] The function $\widehat{\Theta}$ is well-defined and differentiable on $\text{Sym}^2(\mathbb{R}^{n+1} \setminus \mathcal{S})$, and for each $c \in (-(n+1)\pi/2, (n+1)\pi/2)$ such that $c \equiv \theta \mod 2\pi$, the set*

$$\boldsymbol{\mathcal{F}_c} = \{A \in \text{Sym}^2(\mathbb{R}^{n+1}) : \widetilde{\Theta}(A) \geq c\}$$

*is a subequation for the DSL of phase $\theta$.*

**Remark 2.1.2.** *The different choices of $c$ for a given $\theta$ correspond to the different branches of the DSL. The DSL subequation is unique in the sense that it arises as the super-level set of an upper semi-continuous function and not a continuous one. $\widehat{\Theta}$ cannot be extended continuously to all of $\text{Sym}^2(\mathbb{R}^{n+1})$. See [RS, Section 3] for more details.*

The positivity condition defining the space of positive Lagrangians can also be phrased in terms of a subequation, namely the special Lagrangian subequation introduced by Harvey–Lawson,

$$\boldsymbol{F_c} := \{A \in \mathrm{Sym}^2(\mathbb{R}^n) : \mathrm{tr}\tan^{-1}(A) \geq c\}.$$

**Theorem 2.1.3.** *[RS, Corollary 5.6] Let $\theta \in (-\pi, \pi]$, let $D \subset \mathbb{R}^n$ be a domain and let $k \in C^2([0,1] \times D)$. Then $k$ is a solution of the DSL if and only if for each $(t,x) \in [0,1] \times D$,*

$$\nabla^2 k(t,x) \in \boldsymbol{\mathcal{F}_c} \cap -\widetilde{\boldsymbol{\mathcal{F}}}_c = \partial \boldsymbol{\mathcal{F}_c},$$

$$\nabla_x^2 k(t,x) \in \mathrm{int}\left(\boldsymbol{F_{c-\pi/2}} \cap -\boldsymbol{F_{-c-\pi/2}}\right),$$

*for a fixed $c \in (-(n+1)\pi/2, (n+1)\pi/2)$ satisfying $c = \theta + 2\pi k$ with $k \in \mathbb{Z}$.*

## 2.2 The Riemannian DSL subequation

In this section we prove the following.

**Proposition 2.2.1.** *For any Riemannian manifold $M$, the Riemannian manifold $\mathbb{R} \times M$ admits a global Riemannian subequation $\mathcal{F}_c$ locally modelled on the Euclidean degenerate special Lagrangian subequation $\boldsymbol{\mathcal{F}_c}$.*

Proposition 1.1.1 is then an immediate consequence. To prove Proposition 2.2.1 we show that the $(n+1)$-dimensional manifold $\mathbb{R} \times M$ admits a topological $O_n$-structure and that $\boldsymbol{\mathcal{F}_c}$ has compact invariance group $O_n$. This implies that $\boldsymbol{\mathcal{F}_c}$ induces a global equation on $\mathbb{R} \times M$. See Section A.5.

*Proof of Proposition 2.2.1.* To see that $\mathbb{R} \times M$ admits an $O_n$-structure (viewing $O_n \subset O_{n+1}$), observe that because $\mathbb{R} \times M$ is globally a product,

$$T(\mathbb{R} \times M) \cong T\mathbb{R} \oplus TM,$$

and $\mathbb{R}$, being parallelizable, admits a trivial structure. In terms of the metric, since $O_n$-structures are equivalent to Riemannian structures, this represents the fact that $\mathbb{R} \times M$ admits a global product metric.

Now we show that the compact invariance group of $\boldsymbol{\mathcal{F}_c}$ contains $O_n$. Let $S$ denote the elements $A \in \mathrm{Sym}^2(\mathbb{R}^{n+1})$ of the form $A = \mathrm{diag}(0, B)$, for some $B \in \mathrm{Sym}^2(\mathbb{R}^n)$ and set $I_n = \mathrm{diag}(0, I)$. From [RS, Section 3.1] it follows that:

If $A \in S$,
$$\widetilde{\Theta}(A) = \frac{\pi}{2} + \mathrm{tr}\ \arg(I + \sqrt{-1}B).$$

If $A \in \mathrm{Sym}^2(\mathbb{R}^{n+1}) \setminus S$,
$$\widetilde{\Theta}(A) = \sum_{\lambda \in \mathrm{spec}(I_n + \sqrt{-1}A)} m(\lambda)\arg(\lambda).$$

Let $H = \mathrm{diag}(1, h) \in O_{n+1}$, where $h \in O_n$.

When $A \in S$, $A = \mathrm{diag}(0, B)$, so
$$\widetilde{\Theta}(HAH^t) = \frac{\pi}{2} + \mathrm{tr}\ \arg(I + \sqrt{-1}hBh^t)$$
$$= \frac{\pi}{2} + \mathrm{tr}\ \arg(I + \sqrt{-1}B).$$
$$= \widetilde{\Theta}(A).$$

When $A \notin S$, we have
$$I_n + \sqrt{-1}HAH^t = HI_nH^t + \sqrt{-1}HAH^t = H(I_n + \sqrt{-1}A)H^t,$$

and since $O_{n+1} \subset U_{n+1}$ the spectrum of

$$I_n + \sqrt{-1}A \quad \text{and} \quad I_n + \sqrt{-1}hAh^t$$

are the same. Thus, $\widetilde{\Theta}(HAH^t) = \widetilde{\Theta}(A)$. Therefore, the compact invariance group of $\boldsymbol{\mathcal{F}_c}$ contains $O_n \subset O_{n+1}$. □

The special Lagrangian subequation $\boldsymbol{F_c}$ has been studied by Harvey–Lawson. See [HL2, Section 10] and [HL3, Section 14]. Since $\boldsymbol{F_c}$ depends only on the eigenvalues of $A$ it is $O_n$-invariant and caries over to a Riemannian subequation $F_c$ on any $n$-dimensional Riemannian manifold.

## 2.3 The DSL on complex manifolds and higher corank

In unpublished notes [R], Rubinstein showed that the DSL subequation can also be defined in the complex setting, i.e., there is a well defined subequation for the equation

$$\operatorname{Im}\det(I_n + \sqrt{-1}\operatorname{Hess}_{\mathbb{C}}u) = 0, \tag{13}$$

where $\operatorname{Hess}_{\mathbb{C}}$ is the complex $(1,1)$ Hessian. Harvey–Lawson considered the non-degenerate case in [HL3, Section 15].

More specifically, let $k : \mathbb{C}^{n+1} \to \mathbb{C}$ so that

$$(\tau, z_1, \cdots, z_n) \mapsto k(\tau, z_1, \cdots, z_n) = k(t, z_1, \cdots, z_n),$$

where $\tau = t + \sqrt{-1}s$.

Then equation (13) is invariant under the unitary matrices $U_n \subset U_{n+1}$ in the sense that for any $U = \operatorname{diag}(1, V) \in U_{n+1}$, where $V \in U_n$,

$$\operatorname{Im}\det[I_n + \sqrt{-1}U\operatorname{Hess}_{\mathbb{C}}kU^*] = \operatorname{Im}\det[V(I_n + \sqrt{-1}\operatorname{Hess}_{\mathbb{C}}k)V^*] = \operatorname{Im}\det[I_n + \sqrt{-1}\operatorname{Hess}_{\mathbb{C}}k].$$

Any almost complex manifold $X$ admits a topological $U_n$ structure. Since $\mathbb{C}$ is paralellizable, $\mathbb{C} \times X$ also admits a topological $U_n$ structure, viewing $U_n \subset U_{n+1}$. Thus, for any almost complex manifold $X$, there exists a global equation on $\mathbb{C} \times X$ locally modelled on Equation (13). Taking $D \subset X$ to be the domain

$$\mathcal{D} = \{(\tau, z) \in \mathbb{C} \times X : 0 \leq t \leq 1 \text{ and } z \in D\},$$

i.e., an infinite strip of width 1 in the complex plane times $D \subset X$, one can consider the Dirichlet problem for equation (13) on $\mathcal{D}$, with data depending only on the real part.

It was also shown in [R] that there are corresponding subequations for higher co-rank DSL equations on $\mathbb{R}^{n+k}$, with $I_n$ replaced by the diagonal matrix $\operatorname{diag}(0,..,0,1,...1)$ with $k$ zeros and $n$ ones. In an analogous manner, these equations will carry over to equations on $\mathbb{R}^k \times M$, for any $n$-dimensional Riemannian manifold $M$.

# 3 GEOMETRIC MOTIVATION FOR THE RIEMANNIAN DSL

In this section we prove Theorem 1.1.2 and Theorem 1.1.3.

## 3.1 Parallelizable manifolds and almost Calabi–Yau cotangent bundles

Recall that an $n$-dimensional manifold $M$ is *parallelizable* if it admits a global frame field for the tangent bundle. In terms of its topological structure group (see Section A.5), a manifold is parallelizable if it admits an $I$-structure, where $I$ is the trivial subgroup in $\operatorname{GL}(n, \mathbb{R})$. Examples of parallelizable manifolds include all orientable 3-dimensional manifolds and all Lie groups [HL3, Section 5.2].

An *almost Calabi–Yau manifold* is an almost complex Hermitian manifold $X$ with a global section of $\Lambda^{(n,0)}(T^*X)$ whose real part has comass 1 [HL3, Section 1]. This is equivalent to having topological structure group $SU_n$. When $M$ is parallelizable, we can explicitly construct an almost Calabi–Yau structure on $T^*M$ which respects the cotangent bundle fibration. The following construction is based on [HL3, Section 14].

Let $(M, g)$ be a parallelizable Riemannian manifold. Taking a global orthonormal frame $v = (v_1, \ldots, v_n)$, we identify

$$TM = M \times \mathbb{R}^n.$$

Taking the global coframe $w = (w_1, \ldots, w_n)$ to $v$,

$$T^*M = M \times \mathbb{R}^n.$$

Thus, $(v, w)$ forms a global frame for

$$T(T^*M) = T^*M \times \mathbb{R}^{2n},$$

where movement along $M$ is captured by $v$ and movement within the fibre by $w$. Let $v^*$ and $w^*$ denote the dual frames to $v$ and $w$ for the cotangent bundle of $T^*M$.

In terms of this framing, $T^*M$ admits an almost Calabi–Yau structure.

Almost complex structure $J$:
$$Jv_i = w_i, \quad Jw_i = -v_i,$$

Non-vanishing $(n,0)$-form $\Omega$:
$$\Omega = (v_1^* + \sqrt{-1}w_1^*) \wedge \cdots \wedge (v_n^* + \sqrt{-1}w_n^*),$$

Non-degenerate 2-form $\omega$:
$$\omega = \sum_i v_i^* \wedge w_i^*.$$

In general, this structure is not integrable. That is, $J$ is not a (integrable) complex structure and $\omega$ and $\Omega$ are not closed. However, a certain degree of integrability is necessary for Solomon's geometry on the space of positive Lagrangians. For instance, if $J$ not integrable then $\Omega$ will not be closed and the connection on $\mathcal{O}$ may no longer be the Levi–Civita connection. To remedy this, we now consider a special class of parallelizable manifolds on which the above almost Calabi–Yau structure is a true Calabi–Yau structure, as defined in Section 1.2.

A manifold is called *integrably parallelizable* if it admits an atlas of charts such that the differentials of the transition maps are the identity. In terms of topological structure groups, this is equivalent to saying $M$ admits an *integrable I-structure*.

**Theorem 3.1.1.** *[H, Section 1] Let $M$ be connected and parallelizable. Then $M$ is integrably paralellizable if and only if $M$ is open (i.e., non-compact and without boundary) or diffeomorphic to the n-dimensional torus.*

**Example 3.1.2.** *[H] Examples of integrably parallelizable manifolds:*

  i. *Open Lie groups;*

  ii. *Punctured compact connected Lie groups;*

  iii. *Open orientable 3-manifolds;*

  iv. *Diffeomorphic images of the torus;*

  v. *Punctured Stiefel manifolds.*

$S^3$ and $\mathbb{R}P^3$ are examples of parallelizable manifolds that are not integrably parallelizable. We now prove Theorem 1.1.2.

*Proof of Theorem 1.1.2.* We first construct a Calabi–Yau structure on $T^*M$. Since $M$ is integrably parallelizable we have a covering of coordinate charts $\{U_\alpha\}_{\alpha \in A}$ such that the differential of the transition maps is the identity. Let $x$ be coordinates on $U_\alpha$, and consider the induced coordinate charts on $T^*M$:
$$V_\alpha = U_\alpha \times \mathbb{R}^n,$$

where $\mathbb{R}^n$ corresponds to the coordinates $\xi^1, \cdots, \xi^n$ for $dx^1, \cdots, dx^n$. It is easy to see that $T^*M$ is also integrably parallelizable with this corresponding atlas. Let $(v, w)$ to be the global frame for $T(T^*M)$ which on $V_\alpha$ takes the form $v_i = \frac{\partial}{\partial x^i}$ and $w_i = \frac{\partial}{\partial \xi^i}$, and take the almost Calabi–Yau structure on $T^*M$ as defined above. This almost Calabi–Yau structure is integrable: since the complex structure is defined in terms of coordinate vector fields the Nijenhuis tensor vanishes, meaning the complex structure is integrable; it follows that $\Omega$ is holomorphic; and it is clear that $\omega$ is closed and compatible with $J$.

Next, we show that the Riemannian DSL coincides with the geodesic equation for gradient graphs. These computations closely follow [RS, Section 2.4] so we mostly emphasize the differences.

Consider the path of Lagrangians in $T^*M$ given by
$$\Lambda_t = \text{graph}\,(df_t),$$

where $f_t \in C^2(M)$ for $t \in [0, 1]$. Let $g_t : M \to \Lambda_t$, where
$$g_t(p) = (p, df_t|_p).$$

Then, in local coordinates,

$$\frac{dg_t}{dt} = \partial_t df_t = \sum_{i=1}^n \left(\frac{\partial^2 f_t}{\partial t \partial x_i}\right) \frac{\partial}{\partial \xi_i} \quad \text{and} \quad dg_t = \mathrm{I} \oplus \mathrm{Hess}\, f_t\left(\frac{\partial}{\partial x}, \frac{\partial}{\partial x}\right),$$

where $\mathrm{Hess}\, f_t\left(\frac{\partial}{\partial x}, \frac{\partial}{\partial x}\right)$ denotes the Riemannian Hessian in local coordinates on $M$ with respect to the (flat) metric on $M$ induced by the Calabi–Yau structure on $T^*M$. It is obvious that $dg_t$ is the identity in the horizontal direction. To see that it is the Riemannian Hessian in the vertical, we compute the image of $\frac{\partial}{\partial x_i}$ under $dg_t$. Let $\nabla^*$ denote the induced metric connection on the cotangent bundle, given in coordinates by

$$\nabla^*_{\frac{\partial}{\partial x^i}} dx^k = -\Gamma_{ij}^k dx^j.$$

Given $s : (-\epsilon, \epsilon) \to M$ satisfying $s(0) = p$ and $\frac{ds}{dt}(0) = \frac{\partial}{\partial x^i}$, the vertical component at $p \in M$ is given by

$$\frac{d}{dt}\left[\sum_k \frac{\partial f_t}{\partial x^k}(s(t)) dx^k|_{s(t)}\right]_{t=0} = \sum_k \frac{\partial^2 f_t}{\partial x^i \partial x^k}(p) dx^k|_p + \frac{\partial f_t}{\partial x^k}\frac{d}{dt}[dx^k|_{s(t)}]_{t=0}$$

$$= \sum_k \frac{\partial^2 f_t}{\partial x^i \partial x^k}(p) dx^k|_p + \frac{\partial f_t}{\partial x^k}\nabla^*_{\frac{\partial}{\partial x^i}} dx^k|_p$$

$$= \sum_k \frac{\partial^2 f_t}{\partial x^i \partial x^k}(p) dx^k|_p - \frac{\partial f_t}{\partial x^k}\Gamma_{ij}^k dx^j|_p$$

$$= \left(\frac{\partial^2 f_t}{\partial x^i \partial x^j}(p) - \Gamma_{ij}^k \frac{\partial f_t}{\partial x^k}\right) dx^j|_p$$

$$= \text{image under } i^{\text{th}} \text{ row of the matrix } \mathrm{Hess}_p f_t\left(\frac{\partial}{\partial x}, \frac{\partial}{\partial x}\right),$$

where the $(ij)^{\text{th}}$ entry of $\mathrm{Hess}_p f_t\left(\frac{\partial}{\partial x}, \frac{\partial}{\partial x}\right)$ is $\mathrm{Hess}_p f_t\left(\frac{\partial}{\partial x_i}, \frac{\partial}{\partial x_j}\right)$. Note that in this particular construction the $\Gamma_{ij}^k = 0$ as the metric on $M$ is flat. This more general computation will be relevant in the next section where the metric on the base manifold is not flat.

Now expressing $g_t$ with respect to the global frame $(v, w)$,

$$g_t(p) = (p, df_t|_p(v)) = (p, df_t|_p(v_1), \cdots, df_t|_p(v_n)),$$

where $df_t|_p(v_i)$ is the coefficient for $w_i$,

$$\frac{dg_t}{dt} = \partial_t df_t|_p(v) = \sum_{i=1}^n \partial_t df_t(v_i) w_i, \quad \text{and} \quad dg_t = \mathrm{I} \oplus \mathrm{Hess}_p f_t(v, v).$$

Plugging into equation (6),

$$g_t^* \iota_{dg_t/dt} \omega(\cdot) = \omega\left(\frac{dg_t}{dt}, dg_t(\cdot)\right) = -\sum_j w_j^* \otimes v_j^* \left(\sum_{i=1}^n \partial_t df_t(v_i) w_i, dg_t(\cdot)\right) = -d\dot{f}_t(\cdot),$$

where $\dot{f}_t(p)$ denotes the derivative of $f_t(p)$ with respect to $t$, giving us

$$h_t \circ g(t, p) = -\dot{f}_t(p). \tag{14}$$

We then compute

$$g_t^* \iota_{dg_t/dt} \Omega = \sum_{i=1}^n \Omega\left(dg_t/dt, dg_t(v_1 \wedge \cdots \wedge \widehat{v_i} \wedge \cdots \wedge v_n)\right) w_1 \wedge \cdots \wedge \widehat{w_i} \wedge \cdots \wedge w_n.$$

$$= \sum_{i=1}^n \det B_i\, w_1 \wedge \cdots \wedge \widehat{w_i} \wedge \cdots \wedge w_n,$$

where $B_i$, $i = 0, \ldots, n$, is the $n$-by-$n$ matrix obtained by removing the $(i+1)$-th column from the $n$-by-$(n+1)$ matrix

$$B = \left[\sqrt{-1}\partial_t df_t|_x(v) \ \mid \ I + \sqrt{-1}\mathrm{Hess}_x f_t(v, v)\right].$$

Similarly,

$$g_t^*\Omega = \det\left[I + \sqrt{-1}\text{Hess}_x f_t(v,v)\right] w_1 \wedge \cdots \wedge w_n$$
$$= \det B_0 \; w^1 \wedge \cdots \wedge w^n$$

From here, the analysis is the same as that in [RS, Section 2.4]. Solving for the vector field $\zeta_t$ we get

$$\zeta_t = \sum_{i=1}^n a^i(t,p) v_i, \quad \text{where} \quad a^i(t,p) = -(-1)^i \frac{\text{Re}\,(e^{-\sqrt{-1}\theta}\det B_i)}{\text{Re}\,(e^{-\sqrt{-1}\theta}\det B_0)}.$$

Thus, the geodesic equation (9) becomes

$$\text{Im}\, e^{-\sqrt{-1}\theta}\det\left[I_n + \sqrt{-1}\text{Hess}f(\overline{v},\overline{v})\right] = 0, \tag{15}$$

where $\overline{v} = (\frac{\partial}{\partial t}, v)$ is a global frame on $[0,1] \times M$ and the Hessian of $f$ is taken with respect to $t$ and $x$. The positivity condition (3) on the Lagrangians implies that $\text{Re}\,(e^{-\sqrt{-1}\theta}\det B_0) > 0$, or

$$\text{Re}\,\det\left[I + \sqrt{-1}\text{Hess}f_t(v,v)\right] > 0. \tag{16}$$

$\square$

### 3.2 Calabi–Yau torus fibrations

In this section, inspired by a paper of Leung–Yau–Zaslow [LYZ], we consider Calabi–Yau manifolds which admit a smooth torus fibration. That is, a Calabi–Yau manifold $X$ which is actually a fibred manifold $\pi : X \to M$, where for any $p \in M$, $\pi^{-1}(p) = \mathbb{T}^n$. We show that the geodesic equation for positive Lagrangian sections corresponds to the Riemannian DSL on $[0,1] \times M$, proving Theorem 1.1.4. We begin by summarizing the calculations of [LYZ, Section 3].

Let $X$ be a Calabi–Yau $n$-fold admitting a smooth torus fibration over a base manifold $M$, possibly compact. And let $\phi$ be a $\mathbb{T}^n$-invariant Kähler potential on $X$. That is, $\phi(x^j, y^j) = \phi(x^j)$, where $y$ are local coordinates on the fiber and $x$ local coordinates on the base. The coordinates $z^j = x^j + iy^j$ are holomorphic on $X$, and the Kähler metric and form are given, respectively, by

$$h = \frac{\partial^2 \phi}{\partial x^i \partial x^j}(dx^i \otimes dx^j + dy^i \otimes dy^j) \quad \text{and} \quad \omega = \frac{\sqrt{-1}}{2}\frac{\partial^2 \phi}{\partial x^i \partial x^j}(dz^i \wedge d\overline{z}^j).$$

By Calabi [C], $X$ is Ricci-flat and $\Omega = dz^1 \wedge \cdots \wedge dz^n$ is covariant constant if and only if $\phi$ satisfies the real Monge–Ampère equation

$$\det \frac{\partial^2 \phi}{\partial x^i \partial x^j} = c, \tag{17}$$

for some constant $c$. Since $\phi$ satisfies (17), the Calabi–Yau condition

$$\frac{\omega^n}{n!} = c\,(-1)^{n(n-1)/2}\left(\frac{\sqrt{-1}}{2}\right)^n \Omega \wedge \overline{\Omega}$$

is satisfied.

Because of the semi-flatness of $h$, $X$ is locally isometric to the tangent bundle $TM$ with the metric induced by $g = \frac{\partial^2 \phi}{\partial x^i \partial x^j}(dx^i \otimes dx^j)$ on $M$. Moreover, if this metric on $M$ is used to identify its tangent and cotangent bundles, then $\omega$ is the standard symplectic form on the cotangent bundle.

Consider a Lagrangian section $C$ of this fibration, locally written as $y(x)$, in $X$. Using the identification with the cotangent bundle and the fact that $C$ is Lagrangian with respect to $\omega$ if and only if it is closed and hence locally exact, it is shown [LYZ] that locally

$$y^j = \phi^{jk}\frac{\partial f}{\partial x^k},$$

for some function $f$, and further computations show

$$dz^1 \wedge \cdots \wedge dz^n|_C = \det\left(I + \sqrt{-1}g^{-1}\text{Hess}\,f\right) dx^1 \wedge \cdots dx^n. \tag{18}$$

*Proof of Theorem 1.1.3.* Let $C_t : M \to X$ be a smooth path of Lagrangian sections, parametrized by $t \in [0,1]$. Take $g_t : M \to C_t$, where $g_t(p) = (p, C_t(p))$. Then, locally, by the above analysis,

$$\frac{dg_t}{dt} = \sum_{j=1}^n \phi^{jk} \frac{\partial^2 f_t}{\partial t \partial x^k} \frac{\partial}{\partial y^j}.$$

By calculations similar to those in Section 3.1 and [RS, Section 2.4],

$$g_t^* \iota_{\frac{dg_t}{dt}} \omega = -d\dot{f}_t(x), \qquad g_t^* \iota_{\frac{dg_t}{dt}} \Omega = \sum_{i=1}^n \det B_i \ dx^1 \wedge \cdots \wedge \widehat{dx^i} \wedge \cdots \wedge dx^n,$$

and

$$\begin{aligned} g_t^* \Omega &= \det \left[ \mathrm{I} + \sqrt{-1} g^{-1} \mathrm{Hess}\ f_t \right] dx^1 \wedge \cdots \wedge dx^n \\ &= \det B_0\ dx^1 \wedge \cdots \wedge dx^n, \end{aligned}$$

where $B_i$, $i = 0, \ldots, n$, is the $n$-by-$n$ matrix obtained by removing the $(i+1)$-th column from the $n$-by-$(n+1)$ matrix

$$B = \left[ \sqrt{-1} \partial_t df_t \ | \ \mathrm{I} + \sqrt{-1} g^{-1} \mathrm{Hess}\ f_t \right].$$

Taking $\mathfrak{g} = e \oplus g$, where $e$ is the Euclidean metric on $\mathbb{R}$, we can express the positivity and geodesic conditions, respectively, as

$$\mathrm{Re}\ e^{-i\theta} \det \left[ \mathrm{I} + \sqrt{-1} g^{-1} \mathrm{Hess}\ f_t \right] > 0 \quad \text{and} \quad \mathrm{Im}\ e^{-\sqrt{-1}\theta} \det \left[ \mathrm{I}_n + \sqrt{-1} \mathfrak{g}^{-1} \mathrm{Hess}\ f \right] = 0,$$

where the second Hessian is taken with respect to $\mathfrak{g}$ and $f$ is viewed as a function on $[0,1] \times M$. Finally, choosing a local orthonormal admissible frame $v = (v_1, \ldots, v_n)$ on $M$ and extending this to $\mathbb{R} \times M$, by $\overline{v} = (\partial_t, v)$ we rewrite these conditions as

$$\mathrm{Re}\ e^{-\sqrt{-1}\theta} \det \left[ \mathrm{I} + \sqrt{-1} \mathrm{Hess}\ f_t(v,v) \right] > 0 \quad \text{and} \quad \mathrm{Im}\ e^{-\sqrt{-1}\theta} \det \left[ \mathrm{I}_n + \sqrt{-1} \mathrm{Hess}\ f(\overline{v}, \overline{v}) \right] = 0.$$

□

# 4 DIRICHLET DUALITY THEORY WITH WEAK BOUNDARY ASSUMPTIONS

Let $F$ be a Riemannian subequation on a manifold $\mathcal{M}$. In this section we extend Dirichlet duality theory to include certain domains $U \subset \mathcal{M}$ with corners.

Following the conventions of Joyce [J], given a manifold with corners $U$, the boundary $\partial U$ is itself a manifold with corners, equipped with a map

$$i_U : \partial U \to U,$$

which may not be injective. The manifold $U$ is said to be a *manifold with embedded corners* if $\partial U$ can be written as the disjoint union of a finite number of open and closed subsets on each of which $i_U$ is injective. A function $\phi$ on $\partial U$ is called *consistent* if it is constant on the fibres of $i_U$.

In the case that $\mathcal{M} = \mathbb{R}^n$, Rubinstein–Solomon extended Dirichlet duality to include such domains.

**Theorem 4.0.1.** *[RS, Theorem 7.8] Let $F$ be a subequation in $\mathrm{Sym}^2(\mathbb{R}^n)$, and let $U$ be a bounded domain in $\mathbb{R}^n$ such that $\overline{U}$ is a manifold with embedded corners. Let $\phi$ be a consistent continuous function on $\partial U$. Assume $\partial U$ is strictly $(F, \phi)$-convex and strictly $(\widetilde{F}, -\phi)$-convex. Then the $F$-Dirichlet problem for $(U, \phi)$ admits a unique solution in $C^0(\overline{U})$.*

Theorem 4.0.1 was then used to obtain continuous solutions to the Dirichlet problem for the DSL, which, as previously mentioned, is naturally posed on a domain with corners. Our goal here is to achieve a similar extension in the setting of Riemannian manifolds (Theorem 4.2.1) and use it to obtain solutions to the Dirichlet problem for the Riemannian DSL.

We briefly outline our approach and provide context for it in relation to [RS, HL2, HL3]. In Section 4.1, we extend the notion of boundary convexity (in a weakened sense) to domains with corners. This is accomplished by decomposing the boundary into a part that is convex (in the original sense) and a part where given subsolutions are well-behaved. Our definitions come straight from [RS, Section 7.2].

Section 4.2 is then devoted to proving Theorem 4.2.1. Because of the local nature of the arguments used in Dirichlet duality in the Riemannian setting, the proofs in [HL3, Section 12] carry over almost exactly to this setting. This can be contrasted to the Euclidean setting, where the use of global defining functions to construct barriers [HL2, Theorem 5.12] makes this extension more difficult. See [RS, Proposition 7.3].

## 4.1 Weak boundary convexity

Let $U \subset \mathcal{M}$ be a bounded domain, and let $\partial U$ denote the boundary of $\overline{U}$ considered as a manifold with corners.

**Definition 4.1.1.** *The boundary component $\partial U_i$ is called strictly $F$-convex if for each $x \in \partial U_i$, $\partial U_i$ is strictly $F$-convex at $x$ in the sense of Definition A.7.1.*

**Definition 4.1.2.** *Let $\phi \in C^0(\partial U)$ be consistent. A subsolution of the $F$-Dirichlet problem for $(U, \phi)$ is a function $u \in F(U) \cap USC(\overline{U})$ such that $u|_{\partial U} \leq \phi$. A subsolution $u$ for $(U, \phi)$ is called $\delta$-maximal at $p \in \partial U$ if $u(p) \geq \phi(p) - \delta$, and maximal at $p$ if $u(p) = \phi(p)$.*

**Definition 4.1.3.** *We say $\partial U$ is strictly $(F, \phi)$-convex if we can decompose $\partial U$ as the disjoint union $A \cup B$, where $A$ and $B$ are unions of components and satisfy the following:*

1. *For each $p \in A$ and $\delta > 0$ there exists a $C^0(\overline{U})$ subsolution of the $F$-Dirichlet problem for $(U, \phi)$ that is $\delta$-maximal at $p$.*

2. *$B$ is strictly $F$-convex.*

## 4.2 Solution of the Dirichlet problem

The main result of this section is the following extension of Theorem A.8.1 and analogue of Theorem 4.0.1. Here we make the same minor technical assumption that is made for Theorem A.8.1 - that is, we assume $F(U)$ and $\widetilde{F}(U)$ both contain at least one function bounded from below.

**Theorem 4.2.1.** *Suppose $F$ is a subequation on $\mathcal{M}$ for which comparison holds. Let $U \subset \mathcal{M}$ be a bounded domain such that $\overline{U}$ is a manifold with embedded corners, and let $\phi$ be a consistent function on $\partial U$. If $\partial U$ is strictly $(F, \phi)$-convex and strictly $(\widetilde{F}, -\phi)$-convex, then the $F$-Dirichlet problem for $(U, \phi)$ admits a unique solution in $C^0(\overline{U})$.*

**Remark 4.2.2.** *By Theorem A.6.1, the existence of certain subequation-specific $C^2$ functions on $\mathcal{M}$ implies that comparison holds for that subequation. For instance, if $\mathcal{M}$ carries a strictly convex $C^2$ function, then by [HL3, Theorem 9.13] comparison holds for every pure second order subequation on $\mathcal{M}$. Thus, in particular, when $\mathcal{M} = \mathbb{R}^n$ comparison holds for all pure second order subequations.*

The proof of Theorem 4.2.1 is divided into a series of smaller steps, following almost exactly [HL3, Section 12].

**Definition 4.2.3.** *Given a consistent continuous function $\phi$ on $\partial U$, consider the Perron family*

$$\mathscr{F}(\phi) \equiv \{u \in USC(\overline{U}) : u|_U \in F(U) \text{ and } u|_{\partial U} \leq \phi\}$$

*and define the Perron function*

$$u_\phi(x) \equiv \sup\{u(x) : u \in \mathscr{F}(\phi)\}$$

*to be the upper envelope of the Perron family.*

**Proposition 4.2.4 ($F$).** *Let $\phi$ be a consistent continuous function on $\partial U$ and suppose $\partial U$ is strictly $(F, \phi)$-convex at $x_0 \in \partial U$. Then for each $\delta > 0$ small, there exists $w \in \mathscr{F}(\phi)$ such that*

  i. *$w$ is continuous at $x_0$*

  ii. *$w(x_0) \geq \phi(x_0) - \delta$*

  iii. *$w \in F(\overline{U})$.*

**Lemma 4.2.5.** *Let $\phi$ be a consistent continuous function on $\partial U$. Let $x_0$ be a point of a boundary component $\partial U_i \subset \partial U$ that is strictly $F$-convex. Then for each $\delta > 0$ small, there exists $w \in \mathscr{F}(\phi)$ such that*

  i. *$w$ is continuous at $x_0$*

  ii. *$w(x_0) = \phi(x_0) - \delta$*

  iii. *$w \in F(\overline{X})$.*

*Proof.* The proof of this is identical to that of [HL3, Proposition F], as the existence of barriers (Theorem A.7.2) is a purely local condition. □

Clearly, an analogous result holds for strictly $\widetilde{F}$-convex boundary components, providing an element in $\widetilde{\mathscr{F}}(-\phi)$ with the corresponding properties.

*Proof of Proposition 4.2.4.* This follows either by assumption or Lemma 4.2.5. □

**Proposition 4.2.6 ($\widetilde{F}$).** *Let $\phi$ be a consistent continuous function on $\partial U$ and suppose $\partial U$ is strictly $(\widetilde{F}, \phi)$-convex at $x_0 \in \partial U$. Then for each $\delta > 0$ small, there exists $w' \in \widetilde{\mathscr{F}}(-\phi)$ such that*

  i. *$w'$ is continuous at $x_0$*

  ii. *$w'(x_0) \geq -\phi(x_0) - \delta$*

  iii. *$w' \in \widetilde{F}(\overline{U})$.*

*Proof.* Same as Proposition 4.2.4 with an exchange of roles. □

Given a function $f$, let usc $f$ denote its *upper semicontinuous regularization*

$$\text{usc } f := \lim_{\delta \to 0} \sup\{f(y) : y \in U \text{ and } d(x,y) < \delta\},$$

and let lsc $f$ denote its lower semicontinuous regularization, defined analogously.

**Lemma 4.2.7 ($F$).** usc $u_\phi|_U \in F(U)$

*Proof.* The proof of this is identical to that of [HL3, Lemma F]. □

**Lemma 4.2.8 ($\widetilde{F}$).** -lsc $u_\phi|_U \in \widetilde{F}(U)$

*Proof.* The proof of this is identical to that of [HL3, Lemma $\widetilde{F}$]. □

**Corollary 4.2.9 ($F$).** $\phi(x_0) \leq \text{lsc } u_\phi(x_0)$

*Proof.* This follows from Proposition 4.2.4, and is essentially identical to the proof of [HL3, Corollary F]. Since $w \in \mathscr{F}(\phi)$, we have $w \leq u_\phi$ and thus lsc $w \leq$ lsc $u_\phi$. Because $w$ is continuous at $x_0$ and $w(x_0) \geq \phi(x_0) - \delta$,

$$\phi(x_0) - \delta \leq \text{lsc } u_\phi(x_0) \quad \forall \delta > 0 \text{ small}.$$

□

**Corollary 4.2.10 ($\widetilde{F}$).** usc $u_\phi(x_0) \leq \phi(x_0)$

*Proof.* The proof of this is essentially identical to that of [HL3, Corollary $\widetilde{F}$]. Take $u \in \mathscr{F}(\phi)$, arbitrary. Since $w' \leq -\phi$ on $\partial U$, this implies

$$u + w' \leq 0 \quad \text{on } \partial U.$$

Since $w' \in \widetilde{F}(\overline{U})$, by comparison

$$u + w' \leq 0 \quad \text{on } \overline{U}.$$

Thus, $u_\phi + w' \leq 0$ on $\overline{U}$. By the continuity of $w'$ at $x_0$, and the fact that $w'(x_0) \geq -\phi(x_0) - \delta$, we have

$$\text{usc } u_\phi(x_0) \leq -w'(x_0) \leq \phi(x_0) + \delta \quad \forall \delta > 0 \text{ small}.$$

□

From this series of results we can draw the following conclusions.

1. By Corollary 4.2.9 and Corollary 4.2.10, we have lsc $u_\phi = u_\phi =$ usc $u_\phi = \phi$ on $\partial U$. Thus, $u_\phi$ is continuous on $\partial U$.

2. By Corollary 4.2.10 and Lemma 4.2.7, it follows that usc $u_\phi \in \mathscr{F}(\phi)$.

3. And since usc $u_\phi \in \mathscr{F}(\phi)$, this means usc $u_\phi \leq u_\phi$ on $\partial U$. Thus, $u_\phi =$ usc $u_\phi$.

We can now prove Theorem 4.2.1. This proof is identical to that of [HL3, Theorem 12.4].

*Proof of Theorem 4.2.1.* It only remains to show that $u_\phi$ is $F$-harmonic. By Corollary 4.2.9 and Lemma 4.2.8,
$$-\text{lsc } u_\phi \in \widetilde{\mathscr{F}}(-\phi).$$
By conclusion (1) above,
$$-\text{lsc } u_\phi = -u_\phi \quad \text{on } \partial U.$$
Since $u_\phi|_U \in F(U)$, $-\text{lsc } u_\phi|_U \in \widetilde{F}(U)$, and $u_\phi - \text{lsc } u_\phi \leq 0$ on $\partial U$, comparison implies
$$u_\phi - \text{lsc } u_\phi \leq 0 \quad \text{on } \overline{U}.$$
Thus, $\text{lsc } u_\phi = u_\phi$, and so $u_\phi$ is $F$-harmonic. □

# 5 SOLUTION OF THE DIRICHLET PROBLEM FOR THE RIEMANNIAN DSL

In this section we seek unique continuous solutions to the Dirichlet problem for the Riemannian DSL. Our set-up is the following. Let $D \subset M$ be a bounded domain with $\partial D$ smooth and let $\mathcal{D} = (0,1) \times D$, so that $\overline{\mathcal{D}}$ is a compact manifold with embedded corners in $\mathcal{M} = \mathbb{R} \times M$. We assume that both $F_{c-\pi/2}(D)$ and $\widetilde{F}_{c-\pi/2}(D)$ contain at least one $C^2$ function bounded below, where $F_{c-\pi/2}$ is the corresponding special Lagrangian subequation (see Section 3).

**Theorem 5.0.1.** *Suppose comparison holds for the Riemannian DSL subequation $\mathcal{F}_c$ on $\mathcal{M}$ and that $\partial D$ is strictly $F_{c-\frac{\pi}{2}}$, $\widetilde{F}_{c+\frac{\pi}{2}}$ convex. Let $\phi \in C^0(\partial \mathcal{D})$ be consistent and affine in t when restricted to $[0,1] \times \partial D \subset \partial \mathcal{D}$. Consider the following hypotheses:*

1. *$c > -\frac{\pi}{2}$ and for each $i \in \{0,1\}$,*
$$\phi_i := \phi|_{\{i\} \times D} \in C^2(D) \cap F_{c-\frac{\pi}{2}}(D).$$

2. *For each $i \in \{0,1\}$,*
$$\phi_i \in C^2(D) \cap F_{c-\frac{\pi}{2}}(D) \cap -F_{-c-\frac{\pi}{2}}(D).$$

*If either (1) or (2) holds, there exists a unique solution in $C^0(\overline{\mathcal{D}})$ to the $\mathcal{F}_c$-Dirichlet problem for $(\mathcal{D}, \phi)$.*

**Remark 5.0.2.** *a. The boundary assumptions hold for any $D$ such that $\partial D$ is strictly convex, in the sense that all of the eigenvalues of the second fundamental form $II_{\partial D}$ are strictly positive. See [HL3, Proposition 11.4 and Example 14.9].*

*b. If $M$ (and hence $\mathcal{M}$) carries a strictly convex $C^2$ function, then by [HL3, Theorem 9.13] comparison holds for every pure second order subequation on $\mathcal{M}$ and thus for the Riemannian DSL. In particular, if $M$ is a complete simply-connected Riemannian manifold with non-positive sectional curvature, then the square of the distance function from a fixed point is convex and thus comparison holds. Taking hypothesis 2., gives Theorem 1.1.4.*

*c. It may happen that $M$ does not admit a $C^2$ convex function, but that some $\overline{D} \subset M$ do. In this case, [HL3, Theorem 9.13] implies that comparison (and the theorem) holds for the Riemannian DSL on $\overline{D}$.*

We prove Theorem 5.0.1 by showing that $\partial \mathcal{D}$ is appropriately convex and then applying Theorem 4.2.1, following as closely as possible the approach in [RS, Section 8]. In Lemma 5.1.1 and Lemma 5.1.2, we construct subsolutions to the DSL that are maximal on various parts of the boundary. These are the analogues of [RS, Lemma 8.3] and [RS, Lemma 8.4], respectively.

The proof of Lemma 5.1.1 is essentially identical to that of [RS, Lemma 8.3] (note that our initial data is $C^2$). However, the proof of [RS, Lemma 8.4] does not carry over to a proof of Lemma 5.1.2. This is due to the absence of appropriate global defining functions in the Riemannian setting. Instead, Lemma 5.1.2 is proved by combining the techniques used to prove [RS, Lemma 8.4] and [HL3, Proposition F].

Lemma 5.1.3 then uses both Lemma 5.1.1 and Lemma 5.1.2 to to show that $\partial \mathcal{D}$ is both $(\mathcal{F}_c, \phi)$ and $(\widetilde{\mathcal{F}}_c, -\phi)$ strictly convex. We omit its proof as it is identical to [RS, Lemma 8.5].

## 5.1 Proof of Theorem 5.0.1

Given $\phi_i \in C^2(D)$ (as above), define $v_i \in C^0(\overline{\mathcal{D}})$, by

$$v_0 = \phi_0 - Ct, \quad v_1 = \phi_1 - C(1-t), \tag{19}$$

where $t$ is the coordinate on $\mathbb{R}$.

**Lemma 5.1.1.** *Suppose $\phi_i \in C^2(D) \cap F_{c-\frac{\pi}{2}}(D)$. For each $i \in \{0,1\}$, the function $v_i$ is of type $\mathcal{F}_c$.*

*Proof.* Suppose $v_0$ is not of type $\mathcal{F}_c$. Then, by the definition of $\mathcal{F}_c$ (see Section 2.1), there is a point $(t,x) \in \mathcal{D}$ such that

$$\widetilde{\Theta}\left(\operatorname{Hess}_{(t,x)} v_0(\bar{e}, \bar{e})\right) = \widetilde{\Theta}\left(\operatorname{diag}[0, \operatorname{Hess}_x \phi_0(e,e)]\right) < c,$$

where $\bar{e} = (\partial_t, e)$ is an admissible frame near $(t,x)$. Thus, by the definition of $\widetilde{\Theta}$,

$$\operatorname{tr} \tan^{-1} \operatorname{Hess}_x \phi_0(e,e) < c - \pi/2.$$

However, $\phi_0$ is of type $F_{c-\pi/2}$, so this is a contradiction. The same argument holds for $v_1$. □

**Lemma 5.1.2.** *Let $\phi \in C^0(\partial \mathcal{D})$ be consistent and affine in $t$ when restricted to $[0,1] \times \partial D$. Let $\delta > 0$ and let $(t_0, x_0) \in [0,1] \times \partial D$. If $\partial D$ is $F_{c-\pi/2}$ strictly convex, then there exists a subsolution to the $\mathcal{F}_c$ Dirichlet problem for $(\mathcal{D}, \phi)$ that is $\delta$-maximal at $(t_0, x_0)$.*

*Proof.* Since the boundary of $D$ is strictly $F_{c-\pi/2}$ convex at $x_0$, by Theorem A.7.2 there exists a local defining function $\rho$ for $\partial D$ near $x_0$ which defines a barrier for $F_{c-\pi/2}$ at $x_0$. That is, there exists $C_0 > 0$, $\epsilon > 0$, and $r > 0$ such that in local coordinates the functions

$$\beta_i(x) = \phi_i(x_0) - \delta + C\left(\rho(x) - \epsilon \frac{|x - x_0|^2}{2}\right)$$

are strictly $F_{c-\pi/2}$ subharmonic on $B(x_0, r)$ for all $C \geq C_0$. Here we have written $\phi_i$ to mean $\phi|_{\{i\} \times D}$, for $i = 0, 1$.

By the continuity of $\phi$, we can shrink $r > 0$ so that

$$\phi_i(x_0) - \delta < \phi_i(x) \text{ on } \partial D \cap B(x_0, r).$$

Let $\psi \in F_{c-\pi/2}(D)$ be bounded below, and pick $N > \sup_{\partial D} |\phi_i| + \sup_{\overline{D}} \psi$ so that

$$\psi - N < \phi_i - \delta \text{ on } \partial D.$$

Choose $C$ sufficiently large so that on $(B(x_0, r) \setminus B(x_0, r/2)) \cap \overline{D}$

$$\beta_i < \psi - N$$

and on $B(x_0, r/2) \cap \overline{D}$

$$\beta_i < \phi_i(x).$$

Note that since $\rho$ is a boundary defining function it is negative inside $D$, where defined. As $\phi$ is affine in $t$ along the boundary of $\mathcal{D}$, it follows that on $[0,1] \times (B(x_0, r/2) \cap \partial D)$

$$\beta(t,x) = \phi(t, x_0) - \delta + C\left(\rho(x) - \epsilon \frac{|x - x_0|^2}{2}\right) \leq \phi(t,x),$$

and on $\partial \mathcal{D}$

$$(\psi - N)(t,x) = (\psi - N)(x) \leq \phi(t,x).$$

Now set $w(t,x) := \max\{\beta, (\psi - N)\}$. Then, for every $t$, $w(t,x)$ is equal to $\beta(t,x)$ near $x_0$ and equal to $\psi - N$ outside $B(x_0, r/2)$. Since

$$\operatorname{Hess} \beta(t,x) = \operatorname{diag}(0, \operatorname{Hess} \beta(x)) \quad \text{and} \quad \operatorname{Hess}(\psi - N)(t,x) = \operatorname{diag}(0, \operatorname{Hess}(\psi - N)(x)),$$

it follows that $\beta(t,x)$ and $(\psi - N)(t,x)$ are $\mathcal{F}_c$-subharmonic. Thus, $w(t,x)$, the max of two $\mathcal{F}_c$-subharmonic functions, is also of type $\mathcal{F}_c$ by Theorem A.3.1.

Since $w(t,x)$ is equal to $\beta(t,x)$ near $x_0$ it is immediate that $w(t_0, x_0) = \phi(t_0, x_0) - \delta$. □

**Lemma 5.1.3.** *Let $\mathcal{D}$ and $\phi$ be as in Theorem 5.0.1. Then $\partial \mathcal{D}$ is $(\mathcal{F}_c, \phi)$ strictly convex and $(\widetilde{\mathcal{F}}_c, -\phi)$ strictly convex.*

*Proof.* The proof of this is essentially identical to the proof [RS, Lemma 8.5]. □

*Proof of Theorem 5.0.1.* Combine Lemma 5.1.3 and Theorem 4.2.1. □

# A  DIRICHLET DUALITY ON RIEMANNIAN MANIFOLDS

This appendix summarizes the relevant terminology and results of Dirichlet duality [HL3].

## A.1  The second-order jet bundle

Let $X$ be a smooth $n$-dimensional manifold. The *second-order jet bundle* $J^2(X) \to X$ is the bundle whose fibre at a point $x \in X$ is the quotient $J^2_x = C^\infty_x / C^\infty_{x,3}$, where $C^\infty_x$ denotes the germs of smooth functions at $x$ and $C^\infty_{x,3}$ the subspace of germs which vanish to order 3 at $x$.

If $X$ carries a Riemannian metric then the Riemannian Hessian, defined for any $C^2$ function $u$ and vector fields $V$ and $W$ on $X$ by $(\text{Hess } u)(V,W) := V(Wu) - (\nabla_V W)(u)$, is a section of $\text{Sym}^2(T^*X)$. The following is a well-known result concerning the Riemannian Hessian [HL3, Section 4].

**Theorem A.1.1** (The canonical splitting)**.** *The Riemannian Hessian provides a bundle isomorphism*

$$J^2(X) \to \mathbb{R} \oplus T^*X \oplus \text{Sym}^2(T^*X) \quad \text{by mapping} \quad J^2_x u \to (u(x), (du)_x, \text{Hess}_x u)$$

*for a $C^2$ function $u$ at $x$.*

## A.2  Subequations

Let $P = \{A \in \text{Sym}^2(T^*_x X) : A \geq 0\}$. A subset $F \subset J^2(X)$ satisfies the *Positivity Condition* (P) if

$$F + P \subset F.$$

Take the canonical splitting $J^2(X) = \mathbb{R} \oplus J^2_{\text{red}}(X)$, where $\mathbb{R}$ denotes the 2-jets of locally constant functions and $J^2_{\text{red}}(X)_x \equiv \{J^2_x u : u(x) = 0\}$ is the space of *reduced 2-jets* at $x$, and define $N \subset \mathbb{R} \subset J^2(X)$ to have fibres $N_x = \mathbb{R}^- = \{c \in \mathbb{R} : c \leq 0\}$. A subset $F \subset J^2(X)$ satisfies the *Negativity Condition* (N) if

$$F + N \subset F.$$

A subset $F \subset J^2(X)$ satisfies the *Topological Condition* (T) if

$$(i) \ F = \overline{\text{Int } F}, \quad (ii) \ F_x = \overline{\text{Int } F_x}, \quad (iii) \ \text{Int } F_x = (\text{Int } F)_x.$$

The main existence and uniqueness results for Dirichlet duality assume that $F$ satisfies (P), (T), and (N), so this is formalized as follows. A *subequation $F$ on a manifold $X$* is a subset $F \subset J^2(X)$ satisfying conditions (P), (T), and (N).

## A.3  $F$-subharmonic functions

Let $F \subset J^2(X)$ be closed. The function $u \in C^2(X)$ is $F$-*subharmonic* if its 2-jet satisfies $J^2_x u \in F_x$, for all $x \in X$, and *strictly $F$-subharmonic* if its 2-jet satisfies $J^2_x u \in (\text{Int} F)_x$, for all $x \in X$. This definition extends to the larger class of upper semi-continuous functions on $X$ taking values in $[-\infty, \infty)$, $\text{USC}(X)$, in a viscosity-like way: $u \in \text{USC}(X)$ is said to be $F-$subharmonic if for each $x \in X$ and each function $\phi$ which is $C^2$ near $x$, one has that

$$\{u \leq \phi \text{ near } x_0 \text{ and } u(x_0) = \phi(x_0)\} \implies J^2_x \phi \in F_x.$$

The set of all such functions is denoted by $F(X)$.

**Theorem A.3.1** (Remarkable Properties of $F$-Subharmonic Functions)**.** *Let $F$ be an arbitrary subequation.*

*(Maximums) If $u, v \in F(X)$, then $w = \max\{u, v\} \in F(X)$.*

*(Coherence) If $u \in F(X)$ is twice differentiable at $x \in X$, then $D^2_x u \in F_x$.*

*(Decreasing Sequences) If $\{u_j\}$ is decreasing sequence of functions in $F(X)$ then limit is of type $\mathcal{F}$.*

*(Uniform Limits) If $\{u_j\}$ is a sequence of functions in $F(X)$ that converges uniformly on compact sets then the limit is if type $F$.*

*(Families Locally Bounded Above) If $\mathcal{F} \subset F(X)$ is a family which is locally uniformly bounded above. Then the USC regularization $v^*$ of the upper envelope $v(x) = \sup_{f \in \mathcal{F}} f(x)$ belongs to $F(X)$.*

Given a subset $F \subset J^2(X)$ the *Dirichlet dual* $\widetilde{F}$ of $F$ is defined by $\widetilde{F} = \sim (-\text{Int} F) = -(\sim \text{Int} F)$, and a function $u$ is $F$-*harmonic* if $u \in F(X)$ and $-u \in \widetilde{F}(X)$.

## A.4 Local trivialization

When $X = \mathbb{R}^n$ the 2-jet bundle is canonically trivialized by $J_x^2 u = (u(x), D_x u, D_x^2 u)$, where

$$D_x u = \left(\frac{\partial u}{\partial x_1}(x), ..., \frac{\partial u}{\partial x_n}(x)\right) \text{ and } D_x^2 u = \left(\frac{\partial^2 u}{\partial x_i \partial x_j}(x)\right).$$

Thus, for any open subset $X \subset \mathbb{R}^n$ there is a canonical trivialization

$$J^2(X) = X \times \mathbb{R} \times \mathbb{R}^n \times \mathrm{Sym}^2(\mathbb{R}^n), \quad \text{with fibre} \quad \mathbf{J}^2 = \mathbb{R} \times \mathbb{R}^n \times \mathrm{Sym}^2(\mathbb{R}^n).$$

The notation $J = (r, p, A) \in \mathbf{J}^2$ will be used for the coordinates on $\mathbf{J}^2$. Any subset $\mathbf{F} \subset \mathbf{J}^2$ which satisfies conditions (P), (N), and (T) determines a *Euclidean subequation* on any open subset $X \subset \mathbb{R}^n$ by setting $F = X \times \mathbf{F} \subset J^2(X)$. This subequation is often referred to as just $\mathbf{F}$.

Let $e = (e_1, ..., e_n)$ be a choice of local framing of the tangent bundle $TX$ on some neighborhood $U \subset X$. With this framing the canonical splitting determines a trivialization of $J^2(U)$ given at $x \in U$ by

$$\Phi^e : J_x^2(U) \to \mathbb{R} \oplus \mathbb{R}^n \oplus \mathrm{Sym}^2(\mathbb{R}^n), \quad \text{defined by} \quad \Phi^e(J_x^2(u)) \equiv (u, e(u), (\mathrm{Hess}\, u)(e, e)),$$

where $e(u) = (e_1 u, ..., e_n u)$ and $(\mathrm{Hess}\, u)(e, e)$ is the $n \times n$–matrix with entries $(\mathrm{Hess}\, u)(e_i, e_j)$.

## A.5 Riemannian G-manifolds and Riemannian G-subequations

The general linear group $\mathrm{GL}_n(\mathbb{R})$ has a natural action on the fibre $\mathbf{J}^2$ given by

$$h(r, p, A) = (r, hp, hAh^t) \quad \text{for } h \in \mathrm{GL}_n(\mathbb{R}).$$

For each Euclidean subequation $\mathbf{F} \subset \mathbf{J}^2$ this action determines a *compact invariance group*

$$G(\mathbf{F}) = \{h \in O_n : h(\mathbf{F}) = \mathbf{F}\}.$$

Fix a subgroup $G \subset O_n$. A *topological G-structure* on $X$ is a family of smooth local trivializations of $TX$ over open sets in a covering $\{U_\alpha\}$ of $X$ with $G$-valued transition functions. A *Riemannian G-manifold* is a Riemannian manifold equipped with a topological $G$-structure.

**Lemma A.5.1.** *[HL3, Lemma 5.2] Suppose F is a Euclidean subequation with compact invariance group $G$ and $X$ is a Riemannian G-manifold. For $x \in X$, the condition on a 2-jet $J \equiv J_x^2 u$ that*

$$\Phi^e(J) \equiv (u(x), e_x(u), (\mathrm{Hess}\,_x u)(e, e)) \in \mathbf{F}$$

*is independent of the choice of G-frame $e$ at $x$. Hence there is a well-defined subset $F \subset J^2(X)$ given by*

$$J \in F_x \iff \Phi^e(J)(x) \in \mathbf{F}.$$

This subset $F \subset J^2(X)$ is a subequation on $X$ and will be called the *Riemannian G-subequation on $X$ with Euclidean model $\mathbf{F}$*.

Now in local coordinates $x = (x_1, ..., x_n)$ on $X$, the Riemannian Hessian takes the following form

$$(\mathrm{Hess}\, u)\left(\frac{\partial}{\partial x_i}, \frac{\partial}{\partial x_j}\right) = \frac{\partial^2 u}{\partial x_i \partial x_j} - \sum_{k=1}^n \Gamma_{ij}^k(x) \frac{\partial u}{\partial x_k},$$

where $\Gamma_{ij}^k$ denote the Christoffel symbols of the Levi–Civita connection. In shorthand,

$$(\mathrm{Hess}\, u)\left(\frac{\partial}{\partial x}, \frac{\partial}{\partial x}\right) = D^2 u - \Gamma_x(Du).$$

**Proposition A.5.2.** *[HL3, Proposition 5.5] Let $F$ be a Riemannian G-subequation on $X$ with Euclidean model $\mathbf{F}$ on a Riemannian G-manifold $X$. Suppose $x = (x_1, ..., x_m)$ is a local coordinate system on $U$ and that $e_1, ..., e_n$ is an admissible G-frame on $U$. Let $h$ denote the $\mathrm{GL}_n$-valued function on $U$ defined by $e = h\frac{\partial}{\partial x}$. Then a $C^2$-function $u$ is F-subharmonic on $U$ if and only if*

$$(u, hDu, h(D^2 u - \Gamma(Du))h^t) \in \mathbf{F} \text{ on } U.$$

## A.6 Comparison and approximation

There is a comparison and approximation theory for subequations, which addresses when the sum of an $F$-subharmonic function and an $\widetilde{F}$-subharmonic function satisfy the maximum principle and when an arbitrary $F$-subharmonic function can be uniformly approximated with strictly $F$-subharmonic functions. We briefly introduce the relevant terminology and an important result.

*Comparison holds for the subequation $F$ on $X$* if for all compact sets $K \subset X$, whenever

$$u \in F(K) \quad \text{and } v \in \widetilde{F}(K),$$

the Zero Maximum Principle holds for $u + v$ on $K$, that is,

$$u + v \leq 0 \text{ on } \partial K \quad \implies \quad u + v \leq 0 \text{ on } K.$$

*Strict approximation holds for $F$ on $X$* if for each compact set $K \subset X$, each function $u \in F(X)$ can be uniformly approximated by strictly $F$-subharmonic functions on $K$. A function $u \in C^2(X)$ is said to be *strictly $F$-subharmonic* on $X$ if $J_x^2 u \in \text{Int } F$ for all $x \in X$. This notion extends to upper semicontinuous functions (see [HL3, Definition 7.4]). Let $F_{\text{strict}}(X)$ denote the set of all upper semicontinuous strictly $F$-subharmonic functions.

A subset $M \subset J^2(X)$ is a *convex monotonicity cone for $F$* if $M$ is a convex cone with vertex at the origin and $F + M \subset F$.

**Theorem A.6.1.** *[HL3, Theorem 10.3] Suppose $F$ is a Riemannian $G$-subequation on a manifold $X$. If $X$ supports a $C^2$ strictly $M$-subharmonic function, where $M$ is a monotonicity cone for $F$, then comparison holds for $F$ on $X$.*

## A.7 Boundary convexity and barriers

Recall the canonical decomposition

$$J^2(X) = \mathbb{R} \oplus J_{\text{red}}^2(X)$$

with fibre coordinates $J \equiv (r, J_0)$. A subequation of the form $\mathbb{R} \oplus F$ with $F \subset J_{\text{red}}^2(X)$ is referred to as a *reduced subequation* or a *subequation independent of the $r$ variable*.

Given a subequation $F \subset J_{\text{red}}^2(X)$ independent of the $r$-variable, the *asymptotic interior* $\vec{F}$ of $F$ is the set of all $J \in J_{\text{red}}^2(X)$ for which there exists a neighborhood $\mathcal{N}(J)$ in the total space of $J_{\text{red}}^2(X)$ and a number $t_0 > 0$ such that

$$t \cdot \mathcal{N}(J) \subset F \quad \text{for all } t \geq t_0.$$

Let $\Omega$ be a domain in $X$ with smooth boundary $\partial \Omega$. A *defining function for $\partial \Omega$* is a smooth function $\rho$ defined on a neighborhood of $\partial \Omega$ such that

$$\partial \Omega = \{x : \rho(x) = 0\}, \quad d\rho \neq 0 \text{ on } \partial \Omega, \quad \text{and } \rho < 0 \text{ on } \Omega.$$

For $x \in \partial \Omega$, $J_x^2 \rho = \{0\} \times J_{\text{red},x}^2 \rho$, so we use the notation $J_x^2 \rho = J_{\text{red},x}^2 \rho$

Given a reduced subequation $F$ on $X$ with asymptotic interior $\vec{F}$ and $\Omega \subset X$ a smoothly bounded domain, the $\partial \Omega$ is called *strictly $F$-convex at $x \in \Omega$* if there exists a local defining function $\rho$ for $\partial \Omega$ near $x$ such that $J_x^2 \rho \in \vec{F}_x$. If this holds at every point $x \in \Omega$ then boundary $\partial \Omega$ is *strictly $F$-convex*.

For general subequations boundary convexity is defined as follows. Given any subequation $F \subset J^2(X)$ there is a family of reduced subequations $F_\lambda \subset J_{\text{red}}^2(X)$, $\lambda \in \mathbb{R}$ defined by

$$\{\lambda\} \times F_\lambda = F \cap \left\{\{\lambda\} \times J_{\text{red}}^2(X)\right\}.$$

**Definition A.7.1.** *[HL3, Definition 11.10] Given a general subequation $F \subset J^2(X)$ and a domain $\Omega \subset X$ with smooth boundary, we say that $\partial \Omega$ is strictly $F$-convex at a point $x$ if $\partial \Omega$ is strictly $F_\lambda$-convex at $x$ for each $\lambda \in \mathbb{R}$. The boundary $\partial \Omega$ is called globally $F$-convex if it is $F$-convex at every $x \in \partial \Omega$.*

The importance of boundary convexity is that it implies the existence of barrier functions at boundary points.

Let $\Omega \subset X$ be a smooth domain and let $\rho$ be a local defining function for $\partial \Omega$ near $x_0 \in \partial \Omega$. Then given any $\lambda \in \mathbb{R}$, $\rho$ defines a $\lambda$-*barrier for $F$ at $x_0 \in \partial \Omega$* if there exists $C_0 > 0$, $\epsilon > 0$, and $r_0 > 0$ such that the function

$$\beta(x) = \lambda + C\left(\rho(x) - \epsilon \frac{|x - x_0|^2}{2}\right) \tag{20}$$

is strictly $F$-subharmonic on $B(x_0, r_0)$ for all $C \geq C_0$. If $F$ is a reduced subequation, then we say that $\rho$ *defines a barrier for $F$ at $x_0$*, since the same $\rho$ works for all $\lambda \in \mathbb{R}$. The following result [HL3, Theorem 11.12] connects boundary convexity to the existence of barriers.

**Theorem A.7.2** (Existence of Barriers). *Suppose $\Omega \subset X$ is a domain with smooth boundary $\partial\Omega$ which is strictly $F$-convex at $x_0 \in \partial\Omega$. Then for each $\lambda \in \mathbb{R}$ there exists a local defining defining function $\rho$ for $\partial\Omega$ near $x_0$ which defines a $\lambda$-barrier for $F$ at $x_0$.*

## A.8 Solution of the Dirichlet problem

Let $\Omega \subset\subset X$. Then $g : \overline{\Omega} \to \mathbb{R}$ is said to *solve the $F$-Dirichlet problem on $\overline{\Omega}$ for boundary values $\phi$* if:

$$(a) \quad g \in C(\overline{\Omega}), \quad (b) \quad g \text{ is } F \text{ harmonic on } \Omega, \quad (c) \quad g = \phi \text{ on } \partial\Omega.$$

Given $\phi \in C(\partial\Omega)$, define the *Perron family*

$$\mathscr{F}(\phi) \equiv \{u \in \text{UCS}(\overline{\Omega}) : u|_\Omega \in F(\Omega) \text{ and } u|_{\partial\Omega} \leq \phi\}$$

and the *Perron function* $u_\phi(x) \equiv \sup\{u(x) : u \in \mathscr{F}(\phi)\}$. Assuming that both $F_{\text{strict}}(\overline{\Omega})$ and $\widetilde{F}_{\text{strict}}(\overline{\Omega})$ contain at least one function bounded below (this assumption is minor - see [HL3, Section 12]), Harvey–Lawson prove the following.

**Theorem A.8.1.** *[HL3, Theorem 13.3] Assume comparison holds for the subequation $F$ on $X$ and the domain $\Omega \subset\subset X$ has smooth boundary. If $\partial\Omega$ is both $F$ and $\widetilde{F}$ strictly convex, then for each $\phi \in C(\partial\Omega)$ the Perron function $u_\phi$ uniquely solves the Dirichlet problem on $\Omega$ for boundary values $\phi$.*

# Acknowledgements

The author would like to thank A. Bhattacharya, R. Harvey, H-J. Hein, R. Hunter, B. Lawson, J. Solomon, and A. Yuval for helpful conversations, T. Darvas for many valuable discussions, and the referees for their insightful comments. The author is especially grateful to Y. Rubinstein for his guidance and encouragement and for introducing him to the DSL. Some of this work took place at MSRI supported by NSF Grant No. DMS-1440140 during the Spring 2016 semester.

UNIVERSITY OF MARYLAND
mdellato@math.umd.edu